\documentclass[hidelinkstrue]{siamart1116}

\usepackage{amsfonts}
\usepackage{graphicx}
\usepackage{epstopdf}
\usepackage{algorithmic}
\ifpdf
\DeclareGraphicsExtensions{.eps,.pdf,.png,.jpg}
\else
\DeclareGraphicsExtensions{.eps}
\fi

\usepackage{newfloat}
\DeclareFloatingEnvironment[name=Algorithm]{algofloat}

\numberwithin{theorem}{section}
\newcommand{\TheTitle}{Interval Superposition Arithmetic} 
\newcommand{\TheAuthors}{Y.~Zha, M.~E.~Villanueva, and B.~Houska}
\headers{\TheTitle}{\TheAuthors}

\date{Received: date / Accepted: date}

\usepackage{amsopn}

\ifpdf
\hypersetup{
  pdftitle={Interval Superposition Arithmetic},
  pdfauthor={Yanlin~Zha, Mario~Eduardo~Villanueva\and Boris~Houska}
}
\fi

\usepackage{mathptmx}  
\usepackage{amssymb,amsmath,amsxtra,euscript,mathrsfs}
\usepackage{theorem}
 \usepackage{color}
\DeclareFontFamily{OT1}{pzc}{}
\DeclareFontShape{OT1}{pzc}{m}{it}{<-> s * [1.10] pzcmi7t}{}
\DeclareMathAlphabet{\mathpzc}{OT1}{pzc}{m}{it}

\newcommand{\order}[1]{\mathbf{O}\left(#1\right)}
\newcommand{\diam}[1]{\mathrm{diam}\left( #1 \right)}
\newcommand{\imag}{\mathpzc{i}}
\renewcommand{\proof}{\noindent \textit{Proof.} \hspace{0.2cm}}
\newtheorem{remark}{Remark}

\newcommand{\inv}{\mathrm{inv}}

\title{Interval Superposition Arithmetic}

\author{
Yanlin Zha\footnotemark[1]
\and Mario E. Villanueva\footnotemark[1] 
\and Boris Houska\footnotemark[1]$^{\ ,}$\footnotemark[2]
}

\begin{document}

\maketitle

\renewcommand{\thefootnote}{\fnsymbol{footnote}}
\footnotetext[1]{School of Information Science and Technology, ShanghaiTech University, China.}
\footnotetext[2]{Corresponding Author, borish@shanghaitech.edu.cn}
\renewcommand{\thefootnote}{\arabic{footnote}}

\begin{abstract}
This paper presents a novel set-based computing method, called interval 
superposition arithmetic, for enclosing the image set of multivariate factorable
functions on a given domain. In order to construct such enclosures, the proposed
arithmetic operates over interval superposition models which are parameterized 
by a matrix with interval components. Every point in the domain of a factorable
function is then associated with a sequence of components of this matrix and the superposition,
i.e. Minkowski sum, of these elements encloses the image of the function at this
point. Interval superposition arithmetic has a linear runtime complexity with
respect to the number of variables. Besides presenting a detailed theoretical analysis 
of the accuracy and convergence properties of interval superposition arithmetic, 
the paper illustrates its advantages compared to existing set arithmetics
via numerical examples. 
\end{abstract}

\begin{keywords}
Set based computing, interval arithmetics
\end{keywords}

\begin{AMS}
65G30, 65G40
\end{AMS}

\section{Introduction}
Tools for constructing enclosures of the image set of nonlinear functions are needed
for a wide variety of numerical computing algorithms. These include global 
optimization based on complete-search~\cite{Floudas2013,Neumaier2004}, robust 
and semi-infinite optimization~\cite{Floudas2007,Mitsos2008}, as well as 
validated integration algorithms~\cite{Villanueva2015a,Neher2007}.
Here, factorable functions~\cite{McCormick1976} are functions that can be represented
as a finite recursive composition of atom operations from a (finite) library
$$\mathcal L = \{ +,-, *, \inv, \sin, \exp, \log, \ldots \} \; .$$
This library typically includes binary sums, binary products, and a number of univariate
atom functions such as univariate inversion, trigonometric functions, exponential functions, logarithms,
and others.

Existing methods for computing enclosures of factorable functions can be divided into three
categories: traditional interval arithmetics and its variants, arithmetics using other
convex sets such as ellipsoids or zonotopes, as well as non-convex set arithmetics~\cite{Chachuat2015}.
Interval arithmetics is one of the oldest and most basic tools for set-based 
computing~\cite{Moore1966,Moore2009}. Unfortunately, one of the main limitations
of standard interval arithmetics is that the computed interval enclosures are often much
wider than the exact range of the given factorable function.
This overestimation effect is mainly caused by the so called \emph{dependency problem},
which appears when multiple occurrences of the same variable (interval) are taken
independently during the computation of the enclosure. 
On the other hand, an advantage of interval arithmetics is its favorable 
computational complexity: the evaluation of an interval extension of a factorable
function usually takes only $2$ to $4$ times longer than a nominal evaluation~\cite{Moore2009}.

One way to generalize interval arithmetics is to replace intervals (or interval vectors)
with more general computer representable convex sets. For example, \mbox{McCormick}
relaxations propagate convex lower and concave upper bounds rather than standard 
intervals~\cite{McCormick1976,Mitsos2009}. \mbox{McCormick's} arithmetic sometimes yields tighter 
bounds, but it is also slightly more expensive than interval arithmetics~\cite{Mitsos2009}.
Another class of convex set based enclosure tools is the so-called ellipsoidal
calculus~\cite{Kurzhanski1997,Villanueva2015}, where multi-dimensional ellipsoids rather
than interval vectors are used in order to represent the set enclosures.
Because the storage of an $n$-dimensional ellipsoid grows quadratically 
with the number of variables, i.e., $\order{n^2}$, ellipsoidal arithmetics are typically
computationally more demanding than standard interval arithmetics, but often yield
much tighter enclosures, especially in the context of validated integration algorithms~\cite{Houska2015}.
Thus, at  least for particular applications, the higher computational effort associated 
to ellipsoidal computations pays out in terms of the accuracy of the enclosure 
set. Other convex enclosure methods use polyhedral sets, which are in general 
even more expensive to store than ellipsoids. Unlike ellipsoids, polytopes can 
be used to represent convex sets with arbitrary precisions by controlling the 
number of facets. Polyhedral relaxations are popular in the field of global 
optimization and are for example used in the 
software tools~\texttt{BARON}~\cite{Sahinidis1996,Tawarmalani2005} and \texttt{GloMIQO}~\cite{Misener2013}.
Another example for an enclosure algorithm based on polyhedral sets is the so-called 
\emph{affine arithmetic}~\cite{Figueiredo2004}, which is based on zonotopes, a 
particular class of point-symmetric polytopes.

A rather apparent disadvantage of all arithmetics based on convex sets is that they 
can, in the best case, represent the convex hull of the image of a given 
factorable function. Consequently, if the exact image set of a factorable function
is non-convex, the benefit of investing into more accurate convex set 
representations, such as zonotopes or even general polytopes with many facets, 
is limited. One way to overcome this limitation is by working with non-convex 
sets, which, in practice, is often done using polynomials. Interval polynomials
or polynomials with interval remainder terms have been in use since their
development in the 1960s~\cite{Moore1966} and 
1980s~\cite{Eckmann1984,Ratschek1984}. These early works have been the basis for
the popular Taylor model arithmetic, which has been developed by Berz and 
coworkers~\cite{Berz1997,Berz1998,Makino1999}. Nowadays, there exist mature tools,
for example the software~\texttt{MC++}~\cite{Mitsos2009}, implementing Taylor 
model arithmetics with arbitrary order. The favorable convergence properties of 
Taylor models on variable domains with small diameter have been analyzed 
thoroughly~\cite{Bompadre2013}. However, the convergence properties of Taylor 
series on wider domains are often less favorable~\cite{Rajyaguru2016}.

One promising direction towards overcoming this limitation of Taylor models is the
ongoing research on so-called Chebychev models. For functions with one or two 
variables Chebychev models can be constructed by the software \texttt{Chebfun} 
as developed by Trefethen and coworkers~\cite{Battles2004,Trefethen2007,Townsend2013}.
Chebychev models for functions with more than two variables are the focus of 
recent research~\cite{Rajyaguru2016}. While computing bounds on convex sets is 
computationally tractable, finding tight bounds of a multivariate polynomial is 
itself a complex task. Here, one way to compute bounds on such polynomials is to
use linear matrix inequalities~\cite{Lasserre2009}. Heuristics for computing 
range bounders for multivariate polynomials can be found in~\cite{Rokne1995}.

The main contribution of this paper is the development of a novel non-convex 
set arithmetic, called \emph{interval superposition arithmetic}, for enclosing 
the image set of factorable functions on a given interval domain.
The paper starts in Section~\ref{sec::introISA} by introducing
\emph{interval superposition models}, a data structure that can be used to
represent piecewise constant enclosure functions. In contrast to
the above reviewed non-convex set based arithmetics the ongoing developments
do not rely on local approximation methods such as variational analysis,
Taylor expansions, or other polynomial approximation techniques.
Instead, Section~\ref{sec::ISArithmetic} presents algorithms for propagating
interval superposition models through the directed acyclic graph of factorable
functions by exploiting partially separable sub-structures. Moreover, we develop
associated remainder bounds by exploiting globally valid algebraic
properties, such as addition theorems, which can be found in Appendix~\ref{sec::bounds}.
A detailed analysis of the local
convergence properties of the proposed arithmetic as well as
results on its global behavior can be found in
Sections~\ref{sec::localAnalysis1} and~\ref{sec::globalAnalysis1}, respectively.
Section~\ref{sec::examples} presents numerical results based on a prototype
implementation of the proposed interval superposition arithmetic,
written in the programming language \texttt{JULIA}. The numerical case
studies show that the proposed arithmetic often yields more accurate enclosures
of factorable functions than existing interval arithmetics and Taylor model based
arithmetics, at least on wider domains.
Section~\ref{sec::conclusions} concludes the paper.

\subsection*{Notation}
We use the symbol
\[
\mathbb I \; = \; \left\{ \, [a,b] \subseteq R \, \mid \, a,b \in \mathbb R, \, 
a \, \leq \, b \right\}.
\]
to denote the set of real valued compact interval vectors.
The notation 
$c + I = I + c$ with $I = [a,b] \in \mathbb{I}$ and 
$c\in\mathbb{R}$ is used to represent the shifted interval 
$[c + a, c + b]$. Similarly, $c I = I c$ denotes the scaled
interval $[c a, c b]$ if $c\geq 0$ and $[c b, c a]$ if $c < 0$.
All other interval operations are assumed to be evaluated by
a simple application of standard interval arithmetic.
For example, we use the shorthand notation
\begin{eqnarray}
[a,b] + [c,d] &=& [a+b,c+d] \notag \\
\left[ a,b \right] * [c,d] &=& [ \min \{ a c, a d, b c, b d  \}, \max \{ a c, a d, b c, b d  \} ] \, \notag \\
\exp([a,b]) &=& [ \exp(a), \exp(b) ] \; , \; \text{etc..} \; . \notag
\end{eqnarray}
A complete list of these standard interval arithmetic operations can be found
in~\cite{Moore1966}.

\section{Interval Superposition Models}
\label{sec::introISA}
Let $f: X \to \mathbb R$ be a given factorable function and
$X = [ \underline x_1, \overline x_1] \times [ \underline x_2, \overline x_2] 
\times \ldots [ \underline x_n, \overline x_n] \in \mathbb I^{n}$ a given interval
domain. A set valued function $F_{f,X}: X \to \mathbb I$ is called an interval 
valued enclosure function of $f$ on the given domain $X$, if it satisfies
\[
\forall x \in X, \qquad f(x) \in F_{f,X}(x) \; .
\]
In the following, coordinate aligned branching is applied in order to cut
the whole domain into smaller intervals of the form
\begin{align}
\label{eq::Xbranches}
X_i^j = \left[ \, \underline x_i + (j-1) h_i , \, \underline x_i + j h_i \, 
\right] \qquad \text{with}
\qquad h_i = \frac{\overline x_i - \underline x_i}{N}
\end{align}
for all $i \in \{ 1, \ldots, n \}$ and all $j \in \{ 1, \ldots, N \}$, where $N$
is an integer that the user can choose. Here, the intervals 
$[\underline x_i,\overline x_i]$ are all cut into $N$ equidistant intervals for 
simplicity of presentation, although the following algorithms can easily be 
generalized for non-equidistant interval branching and for the case that each 
coordinate is not necessarily subdivided into the same number of intervals.
Next, we introduce the basis functions
\begin{align}
\label{eq::basisFunctions}
\phi_i^j(x) = \left\{
\begin{array}{ll}
1 & \text{if} \; x_i \in X_i^j  \\
0 & \text{otherwise}
\end{array}
\right.
\end{align}
for all $i \in \{ 1, \ldots, n \}$ and all $j \in \{ 1, \ldots, N \}$.
Now, the goal is to develop an arithmetic that computes piecewise 
constant enclosure functions of the form
\begin{equation}
\label{eq::superposition1}
F_{f,X}(x) = \sum_{i=1}^n \sum_{j=1}^N \, A_i^j \phi_i^j(x) \; ,
\end{equation}
where the coefficients $A_i^j \in \mathbb I$ are intervals.
The enclosure function $F_{f,X}$ given by~\eqref{eq::superposition1} is called an
\textit{interval superposition model}. This name is motivated by the fact that 
$F_{f,X}(x)$ is represented as a Minkowski sum of $n$ interval valued functions.
Notice that the complexity of storing an interval superposition model
is $2 n N$, as we need to store the upper and lower bounds of the $n N$ intervals
$A_i^j$. The function $F_{f,X}(x)$ is piecewise constant in $x$ and may take 
different interval values on all of its $N^n$ pieces.

In the following, the index $i$ in~\eqref{eq::superposition1} is called the row
index of the  coefficient matrix
\[
A \; = \; \left(
\begin{array}{ccc}
A_1^1 & \ldots & A_1^N \\
\vdots & \ddots & \vdots \\
A_n^1 & \ldots & A_n^N \\
\end{array}
\right) \;.
\]
Similarly, $j$ is called the column index. This matrix notation is introduced in
order to have a convenient storage format for the interval coefficients.

\begin{remark}
\label{rem::redundancy}
Notice that that there is more than one way to represent the same interval 
superposition model. This is mainly due to the fact that the enclosure set
$F_{f,X}(x)$ remains invariant if we pick two pairwise disjoint row indexes, 
$k_1 \neq k_2$, and a constant $c \in \mathbb R$; add the offset $c$ to all 
intervals in the $k_1$-th row; and subtract $c$ from all intervals in the 
$k_2$-th row, i.e.
$$\forall j \in \{ 1, \ldots, N\}, \qquad A_{k_1}^j \leftarrow A_{k_1}^j + c 
\quad \text{and} \quad A_{k_2}^j \leftarrow A_{k_2}^j - c \; .$$
Such redundancies can be removed using a sparse interval matrix $A$, which 
maintains systematically as many zero interval entries as possible. 
\end{remark}

\subsection{Range Bounders}
\label{sec::rangeBounders}
Bounds for the range of an interval superposition model $F_{f,X}$ can be found 
by computing the global minimum and global maximum of the model, i.e. 
\[
\lambda(A) \; / \; \mu(A) \; = \; \min_{x,y} \; / \; \max_{x,y}  \; y \quad \text{s.t.} \quad \left\{
\begin{array}{l}
y \in F_{f,X}(x) \\
x \in X \; .
\end{array}
\right.
\]
The functions $\lambda$ and $\mu$ are called range bounders. Let us denote the row-wise
upper and lower bounds of a given interval matrix $A$ by 
\[
U \left( A_i \right) = \max_{j \in \{ 1, \ldots, N \}} \, \overline A_i^j  \qquad \text{and} \qquad  L \left( A_i \right) = \min_{j \in \{ 1, \ldots, N \}} \, \underline A_i^j \qquad \text{with} \qquad A_i^j = \left[ \underline A_i^j,  \overline A_i^j \right] \; .
\]
The exact range bounders of $F_{f,X}$ can now be evaluated by using the 
following proposition.

\begin{proposition}
	\label{prop::range}
	An interval superposition model $F_{f,X}$ has range $[\lambda(A),\mu(A)]$, 
with 
	\[
  \lambda(A) = \sum_{i=1}^n L(A_i) \quad \text{and} \quad \mu(A) = \sum_{i=1}^n U(A_i) \; .
	\]
\end{proposition}

\proof
The main idea is to exploit complete separability of $F_{f,X}$, i.e.
\begin{equation*}
\forall x \in X, \qquad F_{f,X}(x) = \sum_{i=1}^n \underbrace{\left[ \sum_{j=1}^N \, A_i^j \phi_i^j(x) \right]}_{\text{depends on $x_i$ only}} \; .
\end{equation*}
The definition of the basis functions $\phi_i^j$ in~\eqref{eq::basisFunctions} implies that $\phi_i^j(x)$ depends on $x_i$ only. In other words, the summands in the above expression can be minimized and maximized separately finding the componentwise extrema $L(A_i)$ and $U(A_i)$, respectively. The sum of these extrema corresponds to the exact range bounder of $F_{f,X}$, as stated by the proposition. \hfill\proofbox
\smallskip

An immediate consequence of the above proposition is that if $F_{f,X}$ is an enclosure function of $f$ on $X$, then upper and lower bounds on the function $f$ on the domain $X$ are given by
\[
\forall x \in X, \qquad \lambda(A) = \sum_{i=1}^n L(A_i) \; \leq \; f(x) \; \leq \; \sum_{i=1}^n U(A_i) = \mu(A) \; .
\]
Notice that the cost of evaluating the functions $U$ and $L$ for one row $A_i$ is of order $\order{N}$. Thus, if $F_{f,X}(x)$ is a given superposition model of $f$, the cost of computing the above upper and lower bounds $\mu(A)$ and $\lambda(A)$ is of order $\order{n N}$, as the functions $U$ and $L$ have to be evaluated for all $n$ rows of the coefficient matrix $A$ and added up.

\section{Interval Superposition Arithmetic}\label{sec::ISArithmetic}
This section deals with the propagation of interval superposition models through
a factorable function whose atom operations belong to a library
$$\mathcal L = \{ +, -, *, \inv, \exp, \sin, \log, \ldots \} \, ,$$
which contains bivariate and univariate operators. With respect to the bivariate operators,
we consider only addition and multiplication, as for two given atom operations $g$ and $h$,
their difference and quotient
\begin{equation*}
h(x)-g(x) = h(x) + (-g(x)) \quad \text{and} 
\quad h(x)/g(x) = h(x)*\inv( g(x) ) \; ,
\end{equation*} 
can be obtained by combining binary addition and binary multiplication with 
univariate mirroring and univariate inversion.

\subsection{Univariate Compositions}
\label{sec::1univ}
Let us consider the interval superposition model
\[
F_{h,X}(x) = \sum_{i=1}^n \sum_{j=1}^N \, A_i^j \phi_i^j(x) \;,
\]
of $h: X \to \mathbb R$ on $X \in \mathbb{I}^{n}$.  
Let $g \in \mathcal{L}$ denote a given univariate atom operation.
The goal of this section is to find an interval 
superposition model of the function $f = g \circ h$,
\[
F_{f,X}(x) = \sum_{i=1}^n \sum_{j=1}^N \, C_i^j \phi_i^j(x)\;.
\]
Here, $g \circ h$ denotes the composition of $g$ and $h$, 
$(g \circ h)(x) = g(h(x))$ for all $x$.
The input of a composition rule of a univariate atom operation $g \in \mathcal L$ are the 
coefficients $A_i^j$ and its output are the coefficients $C_i^j$ such that 
whenever $F_{h,X}(x)$ is an enclosure function of $h$ on $X$, then $F_{f,X}$ is 
an enclosure function of $f=g \circ h$ on $X$. Although the particular 
construction of a valid map from $A$ to $C$ differs for each atom function $g$, 
the main concept for computing $C$ is outlined in Algorithm~1.
Notice that the complexity of this algorithm is of order $\order{nN}$.

\begin{algorithm}
	\caption{Composition rule of interval superposition arithmetic}
	\footnotesize
	\bigskip
	\begin{center}
		\begin{minipage}{1\textwidth}
			\textbf{Input:} Interval valued coefficients $A_i^j \in \mathbb I$ of the input model $F_{h,X}$ and an atom function $g \in \mathcal L$.\\
			
			\textbf{Main Steps:}\\[-0.3cm]
			\begin{quote}
				\begin{enumerate}
					\setlength{\itemsep}{2pt}
					\item Choose for all $i \in \{ 1, \ldots, n \}$ suitable central points $a_i \in \mathbb R$ such that
					$$L(A_i) \, \leq \, a_i \, \leq \, U(A_i) \quad \text{and set} \quad \omega = \sum_{i=1}^n a_i \; .$$
					
					\item Choose a suitable remainder bound $r_g(A) \geq 0$ such that
					\begin{eqnarray}
					\left| \sum_{i=1}^n g( \omega + \delta_i ) - (n-1)g(\omega) - g \left( \omega + \sum_{i=1}^n \delta_i \right) \right| \; \leq \; r_g(A) \notag
					\end{eqnarray}
					for all $\delta \in \mathbb R^n$ with $\forall i \in \{ 1, \ldots, n \}, \; \; L(A_i) \leq a_i + \delta_i \leq U(A_i)$.
					
					\item Compute the interval valued coefficients
					\[
					C_i^j = g \left(\omega - a_i + A_i^j \right) - \frac{n-1}{n} g(\omega) \; .
					\]
					for all $i \in \{ 1, \ldots, n\}$ and all $j \in \{ 1, \ldots, N \}$, where $g \left(\omega - a_i + A_i^j \right)$ is evaluated by using traditional interval arithmetic.
					\item Pick a suitable $k \in \{ 1, \ldots, n \}$ and set $C_k^j \leftarrow C_k^j + r_g(A) \cdot [-1,1]$ for all $j \in \{ 1, \ldots, N \}$.
				\end{enumerate}
			\end{quote}
			\smallskip
			
			\textbf{Output:} The coefficients $C_i^j$ of a interval superposition model $F_{f,X}$ of the function $f = g \circ h$.\\
		\end{minipage}
	\end{center}
\end{algorithm}

\begin{theorem}
	Let $F_{h,X}(x) = \sum_{i=1}^n \sum_{j=1}^N \, A_i^j \phi_i^j(x)$ be 
an interval superposition model of $h$ on $X$. If the interval coefficients 
$C_i^j$ are computed by Algorithm~1, the function
	\[
	F_{f,X}(x) = \sum_{i=1}^n \sum_{j=1}^N \, C_i^j \phi_i^j(x)
	\]
	is an interval superposition model of $f = g \circ h$ on $X$.
\end{theorem}

\proof
Let $x \in X$ be any point in the interval $X$. Since 
$F_{h,X} = \sum_{i=1}^n \sum_{j=1}^N \, A_i^j \phi_i^j(x)$ is an 
interval superposition model of the function $h$, there must exist a sequence
of integers $j_1, j_2, \ldots, j_n \in \{ 1, \ldots, N \}$ and associated points
$y_i \in A_i^{j_i}$ such that $h(x) = \sum_{i=1}^n y_i$.
Next, we define $\delta_i = y_i - a_i$ and recall the definition 
$\omega = \sum_{i=1}^n a_i$ from Step~1 of Algorithm~1. These definitions can be
used to write the function $f(x)$ in the form
\begin{eqnarray}
f(x) &=& g \left( h(x) \right) = g \left( \sum_{i=1}^n y_i \right) = g \left( \omega + \sum_{i=1}^n \delta_i \right) \notag \\
&=& \sum_{i=1}^n \left( g \left(\omega + \delta_i \right) - \frac{n-1}{n} g(\omega) \right) - \underbrace{\left( \sum_{i=1}^n g( \omega + \delta_i ) - (n-1)g(\omega) - g \left( \omega + \sum_{i=1}^n \delta_i \right) \right)}_{\in r_g(A) \cdot [-1,1]} \; . \notag
\end{eqnarray}
As we have $\delta_i \in A_i^{j_i} - a_i$, the inclusion 
$g \left( \omega + \delta_i \right) \in g \left( \omega - a_i + A_i^{j_i} \right)$
holds. Consequently, 
\[
f(x) \in \sum_{i=1}^n \left( g \left(\omega - a_i + A_i^{j_i} \right) - \frac{n-1}{n} g(\omega) \right) + r_g(A) \cdot [-1,1] = \sum_{i=1}^n C_i^{j_i} \; .
\]
This implies that $F_{f,X}(x)$, as stated, 
is an interval superposition model of $f = g \circ h$. \hfill\proofbox
\smallskip

The most important steps of Algorithm~1 are Step~1 and Step~2, where central points and an associated remainder bound $r_g(A)$ have to be constructed. This remainder bound is required to satisfy the inequality
\begin{align}
\label{eq::A1bound}
\left| \sum_{i=1}^n g( \omega + \delta_i ) - (n-1)g(\omega) - g \left( \omega + \sum_{i=1}^n \delta_i \right) \right| \; \leq \; r_g(A) 
\end{align}
for all $\delta \in \mathbb R^n$ with $L(A_i) \leq a_i + \delta_i \leq U(A_i)$ 
for all $i \in \{ 1, \ldots, n \}$. Table~\ref{tab::Bounds} lists such central points and remainder bounds for particular atom operations. The corresponding technical derivations of these remainder bounds can be found in 
Appendix~\ref{sec::bounds}.

\begin{table}[!h]
  \footnotesize
	\centering 
	\renewcommand{\arraystretch}{1.7}
	\begin{tabular}{|c|c| c| c|}
		\hline
		Domain & $g(x)$ & Central points & Remainder bound \\
		\hline
		$\mathbb R$ & $-x$ & $a_i = \frac{\mathrm{U}(A_i) + \mathrm{L}(A_i)}{2}$ &
		$\begin{array} {l}
		r_{g}(A) = 0  
		\end{array}
		$ \\
		\hline
		$\mathbb R$ & $x^2$ & $a_i = \frac{\mathrm{U}(A_i) + \mathrm{L}(A_i)}{2}$ & $\begin{array} {l} r_{g}(A) = \sum_{i=1}^n \left( \sigma - s_i \right) s_i \\ 
\text{with} \quad s_i = \frac{\mathrm{U}(A_i)-\mathrm{L}(A_i)}{2} \quad \text{and} \quad \sigma = \sum_{i=1}^n s_i \end{array}$ \\
		\hline
		$\mathbb R_{++}$ & $x^{-1}$ & $a_i = \frac{ \mathrm{L}\left( A_i \right) \mu(A)}{ \lambda(A) + \mu(A) } + \frac{ \mathrm{U}\left( A_i \right) \lambda(A)}{ \lambda(A) + \mu(A) }$ &
		$\begin{array} {l} \\[-0.3cm] r_{g}(A) = \frac{\sum_{i=1}^n s_i \left( \mu(A) - \omega - (U(A_i)-a_i) \right)}{\omega \lambda(A)} \\
		\text{with} \quad s_i = \max \left\{ \, \frac{a_i - L(A_i)}{ \omega - a_i + L(A_i) } \, , \, \frac{U(A_i) - a_i}{ \omega - a_i + U(A_i)} \, \right\}
		\end{array}$ \\
		\hline
		$\mathbb R$ & $e^x$ & $a_i = \log \left( \frac{e^{\mathrm{U}\left( A_i \right)} + e^{\mathrm{L}\left( A_i \right)}}{2} \right)$ & $\begin{array} {l} r_{g}(A) = e^\omega \left[ \prod_{i=1}^n (1+s_i) - \sum_{i=1}^n s_i -1 \right] \\ \text{with} \quad s_i = \frac{e^{\mathrm{U}(A_i)}-e^{\mathrm{L}(A_i)}}{e^{\mathrm{L}(A_i)} + e^{\mathrm{U}(A_i)} } \end{array}$ \\
		\hline
		$\mathbb R_{++}$ & $\log(x)$ & $a_i = \frac{\mathrm{U}(A_i) + \mathrm{L}(A_i)}{2}$ &
		$\begin{array} {l} \\[-0.3cm] r_{g}(A) = - \log \left( 1 - \frac{\prod_{i=1}^n \left( \omega+s_i \right) - \omega^{n-1} \left( \omega + \sum_{i=1}^n s_{i} \right) }{\omega^{n-1} \lambda(A)} \right) \\
		\text{with} \quad s_i = \frac{\mathrm{U}(A_i) - \mathrm{L}(A_i)}{2} \,
		\end{array} \hspace{-0.2cm}$ \\
		\hline
		$\mathbb R$ & $\sin(x)$
		& $a_i = \frac{\mathrm{U}(A_i) + \mathrm{L}(A_i)}{2}$ &
		$\begin{array} {l} \\[-0.3cm] r_{g}(A) = \Omega \left( \prod_{k=1}^n (1+s_k) - \sum_{k=1}^n s_k - 1 \right) \\
		\text{with} \quad \Omega  = |\sin(\omega)| + |\cos(\omega)| \; , \\
		s_i = 2 |\sin([-r_i,r_i])| \; , \quad \text{and} \quad r_i = \frac{\mathrm{U}(A_i) - \mathrm{L}(A_i)}{4} \,
		\end{array} \hspace{-0.2cm}$ \\
		\hline
		$\mathbb R$ & $\cos(x)$
		& same as for $\; \sin(x)$ &
		same as for $\; \sin(x)$ \\
		\hline
		$\left( -\frac{\pi}{2}, \frac{\pi}{2} \right)$ & $\tan(x)$
		& $a_i = \frac{\mathrm{U}(A_i) + \mathrm{L}(A_i)}{2}$ &
		$\begin{array} {l} \\[-0.3cm] 
		r_g(A) = \left| \sum_{i=1}^{n-1} \tan( S_{i+1} ) \tan \left( \sum_{k=1}^i S_k \right)\tan \left( \sum_{k=1}^{i+1} S_k \right) \right. \\
		\; * \left[1 + \tan(\omega) \tan \left( \omega + \Sigma \right) \right] \\
		\; + \sum_{i=1}^n \tan(\omega) \tan \left( S_i \right)  \tan \left( T_i \right) \\
		\; \left. * \left[ 1 + \tan( \omega+S_i)\tan \left( T_i \right) \tan \left( \omega + \Sigma \right) \right] \right|
		\\
		\text{with} \; \; s_i = \frac{U(A_i)-L(A_i)}{2} \; , \; S_i = [-s_{i},s_{i}] \; , \; \sigma = \sum_{i=1}^n s_i \\
		\text{and} \; \; \Sigma = [-\sigma, \sigma]  \; , \; T_i = [-\sigma + s_i, \sigma-s_i]
		\end{array} \hspace{-0.3cm}$ \\
		\hline
	\end{tabular}
	\bigskip
	\begin{center}
	\caption{Central points and remainder bounds for common univariate atom functions.}
	\end{center}
	\label{tab::Bounds}
\end{table}

\begin{remark}
	As discussed in Remark~\ref{rem::redundancy} the proposed interval superposition model storage scheme is redundant with respect to offsets. Consequently, in Step~4 of Algorithm~1 the remainder can in principle be added to any row of the matrix $C$. One possible implementation heuristic is to add the remainder to a row, which contains the intervals with the maximum average diameter. 
\end{remark}

\begin{remark}
Notice that the left column of Table~\ref{tab::Bounds} specifies a domain on which the remainder bound is valid. In some cases this domain can be extended by combining univariate atom operations. For example, an implementation of the function \mbox{$g(x) = x^{-1}$} for negative $x$ is is obtained by combining the atom operations \mbox{$g(x) = x^{-1}$} and \mbox{$g(x) = -x$}. Similarly, the cotangent function can be written in the form $\cot(x) = \tan( \frac{\pi}{2}-x)$. Other functions such as $\sqrt{x} = \exp(0.5 * \log(x))$ can be composed by combining the atom operations in Table~\ref{tab::Bounds}.
\end{remark}

\subsection{Bivariate Compositions}
\label{sec::bivariateComposition}
This section discusses how to construct arithmetic rules for interval 
superpositions for bivariate operators. The addition of two given interval 
superposition models is straightforward. Consider the interval superposition 
models
$$F_{h,X}(x) = \sum_{i=1}^n \sum_{j=1}^N \, A_i^j \phi_i^j(x) \quad \text{and} \quad F_{g,X}(x) = \sum_{i=1}^n \sum_{j=1}^N \, B_i^j \phi_i^j(x)$$
of the given functions $h,g: \mathbb{R}^{n} \to \mathbb R$, on $X\in\mathbb{I}^{n}$. 
Then
$$F_{f,X}(x) = \sum_{i=1}^n \sum_{j=1}^N \, C_i^j \phi_i^j(x) \quad \text{with} \quad C_i^j = A_i^j + B_i^j$$
is an enclosure of the function $f(x) = h(x) + g(x)$.
Algorithm~2 provides a mean to construct an interval superposition model of
$f = g*h$ on $X$, given interval superposition models $F_{g,X}$ and $F_{h,X}$.

\begin{algorithm}
        \caption{Product rule of interval superposition arithmetic}
	\footnotesize
	\bigskip
		\begin{minipage}{1\textwidth}
			\textbf{Input:} Interval valued coefficients $A_i^j \in \mathbb I$ and $B_i^j \in \mathbb I$ of the factors.\\
			
			\textbf{Main Steps:}\\[-0.3cm]
			\begin{quote}
				\begin{enumerate}
					\setlength{\itemsep}{2pt}
					\item Compute the central points
					$$\forall i \in \{ 1, \ldots, n \}, \qquad  a_i = \frac{\mathrm{U}\left( A_i \right) + \mathrm{L}\left( A_i \right)}{2} \quad \text{and} \quad b_i = \frac{\mathrm{U}\left( B_i \right) + \mathrm{L}\left( B_i \right)}{2}$$
					and set
					$$\alpha = \sum_{i=1}^n a_i \; , \; \; \beta = \sum_{i=1}^n b_i \; , \quad \text{and} \quad \gamma = \sum_{i=1}^n a_i b_i \quad \text{as well as} \quad  \omega = \frac{1}{n} \left[ \alpha \beta - \gamma \right] \; .$$
					
					\item Compute the row-wise radii
					\[
					\rho_i(A) = \frac{U(A_i)-L(A_i)}{2} \qquad \text{and} \qquad \rho_i(B) = \frac{U(B_i)-L(B_i)}{2}
					\]
					for all $i \in \{ 1, \ldots, n \}$ as well as the associated remainder bound
					\[
					R(A,B) = \left( \sum_{i=1}^n \rho_i(A) \right) \left( \sum_{i=1}^n \rho_i(B) \right) - \sum_{i=1}^n \rho_i(A)\rho_i(B) \; .
					\]
					
					\item Compute the output coefficients
					\[
					C_{i}^j = \left( A_{i}^j + \alpha - a_i \right) \left( B_{i}^j + \beta - b_i \right) - \left( \alpha - a_i \right)\left( \beta - b_i \right) - \omega
					\]
					for all $i \in \{ 1, \ldots, n\}$ and all $j \in \{ 1, \ldots, N \}$.
					\item Pick a suitable $k \in \{ 1, \ldots, n \}$ and set $C_k^j \leftarrow C_k^j + R(A,B) \cdot [-1,1] $ for all $j \in \{ 1, \ldots, N \}$.
				\end{enumerate}
			\end{quote}
			\smallskip
			
			\textbf{Output:} The coefficients $C_i^j$ of a interval superposition model that encloses the product of the input models.\\
		\end{minipage}
\end{algorithm}

Similar to Algorithm~1, the complexity of Algorithm~2 is of order $\order{n N}$. The validity of the bounds from Algorithm~2 is established in the following theorem.

\begin{theorem}
	Let $F_{h,X}(x) = \sum_{i=1}^n \sum_{j=1}^N \, A_i^j \phi_i^j(x)$ and $F_{g,X}(x) = \sum_{i=1}^n \sum_{j=1}^N \, B_i^j \phi_i^j(x)$ be interval superposition models of $h,g: \mathbb{R}^{n} \to \mathbb R$ on $X\in\mathbb{I}^{n}$. If the coefficients $C_i^j$ are computed by 
Algorithm~2, the function given by 
	\[
	F_{f,X}(x) = \sum_{i=1}^n \sum_{j=1}^N \, C_i^j \phi_i^j(x)
	\]
	is an interval superposition model of the function $f = h * g$ on $X$.
\end{theorem}

\proof
Let $x$ be any point in $X$. Since $F_{h,X}$ and $F_{g,X}$ are interval 
superposition models of the functions $h$ and $g$, there must exist a 
sequence of integers $j_1, j_2, \ldots, j_n \in \{ 1, \ldots, N \}$ and 
associated points $y_i \in A_i^{j_i}$ as well as $z_i \in B_i^{j_i}$ such that
\[
h(x) = \sum_{i=1}^n y_i \quad \text{and} \quad  g(x) = \sum_{i=1}^n z_i \; .
\]
Thus, we have
\begin{align*}
f(x) = h(x)*g(x) \; = \; \left( \sum_{i=1}^n y_i \right)*\left( \sum_{i=1}^n z_i \right) 
=&\hphantom{{}+{}} \sum_{i=1}^n \left[ y_i z_i + y_i \left( \beta - b_i \right) 
+ \left( \alpha - a_i \right) z_i - \omega \right] \\
&- \left( \sum_{i=1}^n \left( y_i - a_i \right) \right) 
\left( \sum_{i=1}^n \left( z_i - b_i \right) \right) \\
&+ \sum_{i=1}^n \left( y_i - a_i \right) \left( z_i - b_i \right) \; .
\end{align*}
Here, the latter equation follows from the addition theorem for the product rule
with 
$$\alpha = \sum_{i=1}^n a_i \; , \; \; \beta = \sum_{i=1}^n b_i \; , \quad \text{and} \quad  \omega = \frac{1}{n} \left[ \left( \sum_{i=1}^n a_i \right) \left( \sum_{i=1}^n b_i \right) - \sum_{i=1}^n a_i b_i \right]$$
The construction of the remainder bound $R(A,B)$ in Step~2 of Algorithm~2 is such that
\[
\left| \left( \sum_{i=1}^n \left( y_i - a_i \right) \right) \left( \sum_{i=1}^n \left( z_i - b_i \right) \right) - \sum_{i=1}^n \left( y_i - a_i \right) \left( z_i - b_i \right) \right| \; \leq \; R(A,B)
\]
for all $y_i \in A_i^{j_i}$ and all $z_i \in B_i^{j_i}$. Consequently, 
\[
f(x) \in \sum_{i=1}^n \left( A_i^{j_i} B_i^{j_i} + A_i^{j_i} \left( \beta - b_i \right) + \left( \alpha - a_i \right) B_i^{j_i} - \omega \right) + R(A,B) \cdot [-1,1] =  \sum_{i=1}^n C_i^{j_i} \; .
\]
This implies that $F_{f,X}(x)$, as stated, is an interval superposition model of 
$f = h*g$.\hfill\proofbox
\smallskip

\subsection{Initialization}
\label{sec::initialization}
Algorithm~1 and~2 can be combined in order to implement the proposed interval 
superposition arithmetic by either operator overloading or source code 
transformation. This is in complete analogy to the implementation of other existing set 
propagation methods operating on the directed acyclic computational graph of the 
given factorable function. The corresponding procedure is initialized by 
constructing (trivial) interval superposition models of all input variables 
$x_i$. As $x_i$ does not depend on other variables its associated interval 
coefficients $A_{k}^j = 0$ can be set to $0$ for all $k \neq i$ and all 
$j \in \{ 1, \ldots, N \}$. The remaining $i$-th row of the interval coefficient
 matrix is initialized by
\[
\forall j \in \{ 1, \ldots, N \}, \qquad A_i^j = X_i^j \; ,
\]
recalling that the branches $X_i^j$ have been defined in~\eqref{eq::Xbranches}.

\section{Properties of interval superposition arithmetic}
This section analyzes the mathematical properties of interval superposition arithmetic. Here, we first analyze the local properties of this arithmetic for small domains $X$. Moreover, Section~\ref{sec::globalAnalysis1} analyzes the global properties and conservatism of the method on large domains.

\subsection{Local overestimation error}
\label{sec::localAnalysis1}
The proposed interval superposition arithmetic is affected by two sources
of overestimation. The first source of overestimation comes from the fact that 
scalar functions, such as $f(x) = x$, can be represented by interval
superposition models with finite accuracy only. However, for Lipschitz continuous
functions, this error is of order $\order{\frac{1}{N}}$ and can be controlled
 by choosing $N$ sufficiently large. Therefore, the focus of the following 
analysis is on the second source of overestimation, namely the remainder bounds 
$r_g(A)$ and $R(A,B)$, needed in Algorithms 1 and 2 respectively.
The following lemma analyzes the local properties of the term that must be bounded by $r_g(A)$.

\smallskip
\begin{lemma}
	\label{lem::RemainderConvergence}
	If the function $g: \mathbb R \to \mathbb R$ is twice continuously differentiable, then
	\begin{eqnarray}
	\left| \sum_{i=1}^n g( \omega + \delta_i ) - (n-1)g(\omega) - g \left( \omega + \sum_{i=1}^n \delta_i \right) \right| \; \leq \; \order{\left[ \mu(A) - \lambda(A) \right]^2} \; . \notag
	\end{eqnarray}
	for all $\delta \in \mathbb R^n$ with $\left| \delta_i \right| \leq U(A_i) - L(A_i)$.
\end{lemma}

\proof
Let $g'$ denote the derivative of the function $g$. As $g$ is twice continuously
differentiable, we can substitute the Taylor expansions
\[
\sum_{i=1}^n g( \omega + \delta_i ) = n g(\omega) + g'(\omega) \sum_{i=1}^n \delta_i + \order{\sum_{i=1}^n  \delta_i^2}
\]
as well as
\[
g\left( \omega + \sum_{i=1}^n \delta_i \right) = g(\omega) + g'(\omega) \sum_{i=1}^n \delta_i + \order{\left[ \sum_{i=1}^n \delta_i \right]^2} \; .
\]
Consequently, we have
\begin{align}
\label{eq::AUX1}
\left| \sum_{i=1}^n g( \omega + \delta_i ) - (n-1)g(\omega) - g \left( \omega + \sum_{i=1}^n \delta_i \right) \right| \; \leq \; \order{\left[ \sum_{i=1}^n \delta_i \right]^2} \; .
\end{align}
We use $\left| \delta_i \right| \leq U(A_i) - L(A_i)$ 
together with the triangle inequality and Proposition~\ref{prop::range} to find
\begin{align}
\label{eq::AUX2}
\left| \sum_{i=1}^n \delta_i \right| \; \leq \; \sum_{i=1}^n \left| \delta_i \right| \; \leq \; \sum_{i=1}^n \left( U(A_i) - L(A_i) \right) \; = \; \mu(A)-\lambda(A) \; .
\end{align}
The statement of the lemma follows now by combining the inequalities~\eqref{eq::AUX1} and~\eqref{eq::AUX2}. \hfill\proofbox

\smallskip
Motivated by Lemma~\ref{lem::RemainderConvergence}, a reasonable requirement on 
$r_g: \mathbb I^{n \times N} \to \mathbb R$ is that it satisfies
\begin{align}
\label{eq::rconv}
\forall A \subseteq \overline D, \qquad r_g(A) \; \leq \; \order{\left[ \mu(A) - \lambda(A) \right]^2} \; ,
\end{align}
where $\overline D \subseteq D$ is a compact subset of the (open) domain $D$ of the atom function $g$. This requirement is satisfied all remainder bounds listed in 
Table~\ref{tab::Bounds} (see Appendix~\ref{sec::bounds} for the details).

\begin{lemma}
	\label{lem::RemainderConvergence2}
	The remainder term $R(A,B)$ of Algorithm~1 satisfies
	\[
	R(A,B) \leq \frac{1}{4} ( \mu(A) - \lambda(A))( \mu(B) - \lambda(B) ) \; .
	\]
\end{lemma}

\proof
The definition of $R(A,B)$ in Step~2 of Algorithm~2 is such that the inequality
\begin{align}
R(A,B) &\leq \left( \sum_{i=1}^n \rho_i(A) \right) \left( \sum_{i=1}^n \rho_i(B) \right) = \frac{1}{4} ( \mu(A) - \lambda(A))( \mu(B) - \lambda(B) )
\end{align}
holds, as stated by the lemma.\hfill\proofbox
\smallskip

The local convergence of interval superposition arithmetic is summarized next.

\begin{theorem}
	\label{thm::localConv}
	Let all atom operations $g \in \mathcal L$ be twice continuously differentiable
and let the remainder bounds $r_g$ of all univariate atom operations 
satisfy~\eqref{eq::rconv}. The maximum distance between the upper and lower bound
of an interval superposition model $F_{f,X}$ computed by the above outlined 
arithmetic rules satisfies
	\[
	\max_{x \in X} \, \left\{  \diam{F_{f,X}(x)} \right\} \; \leq \; \order{ \frac{ \diam{X} }{N} + \left[ \diam{X} \right]^2 } \; ,
	\]
	for all intervals $X \subseteq \overline D$, where $\overline D \subset D$ is a compact subset of an open domain $D$ on which the function $f$ has no singularities. 
\end{theorem}

\proof
The statement of this theorem follows from the fact that variables can be 
represented with accuracy $\order{\frac{ \mathrm{diam}(X) }{N}}$ (induction 
start) while the remainder bound contributions from each atom operation can be 
bounded by expressions of order 
$\order{ \frac{ \mathrm{diam}(X) }{N} + \left[ \mathrm{diam}(X) \right]^2 }$ by 
using the results from Lemma~\ref{lem::RemainderConvergence} 
and~\ref{lem::RemainderConvergence2} (induction step). The details of this 
induction argument are straightforward and skipped for the sake of brevity.
\hfill\hfill\proofbox
\smallskip

At this point, one might argue that the convergence rate of interval superposition
arithmetic is only linear with respect to the diameter of $X$. However, first of
all, the constant in front of the linear term scales with $\frac{1}{N}$ and can 
thus be made arbitrarily small by choosing a sufficiently large $N$. And secondly,
one possible path towards generalizing the above superposition arithmetic could 
be to construct a superposition of Taylor models or other sets rather than 
intervals, if the goal is to move towards better local properties. However, the 
focus of the proposed arithmetic is \textit{not} on the local but rather global
properties of the arithmetic.

\subsection{Global Properties of Interval Superposition Arithmetic}
\label{sec::globalAnalysis1}
In order to discuss the global properties of the arithmetic, we introduce 
the following definition of separability of an interval superposition model.

\begin{definition}
\label{def::separableModel}
An interval superposition model $F_{f,X}(x) = \sum_{i=1}^n \sum_{j=1}^N \, A_i^j \phi_i^j(x)$ is separable, if there exist an integer $k \in \{ 1, \ldots, n \}$
such that $$L(A_i) = U(A_i) \quad \text{for all} \quad i \in \{ 1, \ldots, n \} \setminus \{ k \} \; . $$ 
\end{definition}

An immediate consequence of the initialization routine from 
Section~\ref{sec::initialization} is that the interval superposition model of 
every variable has degree $1$. For the univariate composition rule the following
result can be established.

\begin{lemma}
	\label{lem::zeroRemainder}
	Let the interval superposition model $F_{h,X}$, with interval coefficient $A$, 
of the inner function $h(x)$ in the composition rule (Algorithm~1) be separable. Then
	\[
	\left| \sum_{i=1}^n g( \omega + \delta_i ) - (n-1)g(\omega) - g \left( \omega + \sum_{i=1}^n \delta_i \right) \right| \; = \; 0
	\]
	for all $\delta \in \mathbb R^{n}$ with $L(A_i) \leq a_i + \delta_i \leq U(A_i)$ for all $i \in \{ 1, \ldots, n \}$.
\end{lemma}

\proof
Since $F_{h,X}$ is a separable interval superposition model, we must choose 
$L(A_i) = a_i = U(A_i)$ for all indices $i \in \{ 1, \ldots, n \} \setminus \{ k \}$
for a fixed $k \in \{ 1, \ldots, n \}$. Thus, $\delta_i = 0$ is the only possible
choice for all $i \neq k$. A direct substitution yields
\begin{eqnarray}
& & \left| \sum_{i=1}^n g( \omega + \delta_i ) - (n-1)g(\omega) - g \left( \omega + \sum_{i=1}^n \delta_i \right) \right| \notag \\[0.1cm]
&=& \left| \underbrace{\sum_{i \neq k} g( \omega ) - (n-1)g(\omega)}_{=0} + \underbrace{g( \omega + \delta_k )  - g \left( \omega + \delta_k \right)}_{=0} \right| \, = \, 0  \; , \notag
\end{eqnarray}
which corresponds to the statement of the lemma. \hfill\proofbox
\smallskip

The above lemma implies that the remainder bound function
$r_g: \mathbb I^{n \times N} \to \mathbb R$ can be constructed such that 
$r_g(A) = 0$ whenever the input model 
$F_{f,X}(x) = \sum_{i=1}^n \sum_{j=1}^N \, A_i^j \phi_i^j(x)$ is separable. It
can be checked easily that all remainder bounds from Table~\ref{tab::Bounds} 
have this property.

\begin{lemma} 
	\label{lem::zeroRemainder2}
	If the input models of the product rule from Algorithm~2 are separable with respect to the same index $k$, i.e., if there exists an integer $k \in \{ 1, \ldots, n \}$ such that
	\begin{eqnarray}
	\label{eq::jointlySeparable}
	\forall i \in \{ 1, \ldots, n \} \setminus \{ k \}, \qquad L(A_i) = U(A_i) \quad \text{and} \quad  L(B_i) = U(B_i) \; ,
	\end{eqnarray}
	then the remainder term $R_j(A,B)$ of Algorithm~2 satisfies $R(A,B) = 0$.
\end{lemma}

\proof
If the input models satisfy condition~\eqref{eq::jointlySeparable}, then the equation
\begin{equation*}
\forall i \in \{1, \ldots, n \} \setminus \{ k \} , \qquad  \rho_i(A) = \rho_i(B) = 0
\end{equation*}
is satisfied. A substitution of this equation in the definition of $R$ from Step~2 of Algorithm~2 yields
\begin{eqnarray}
R(A,B) &=& \left( \sum_{i=1}^n \rho_i(A) \right) \left( \sum_{i=1}^n \rho_i(B) \right) - \sum_{i=1}^n \rho_i(A)\rho_i(B) \notag \\[0.16cm]
&=& \rho_k(A) * \rho_k(B) - \rho_k(A) * \rho_k(B) = 0 \; . \notag
\end{eqnarray}
This is the statement of the lemma.\hfill \proofbox
\smallskip

A combination of the above lemmata yields the following global statement about the accuracy of the proposed interval superposition arithmetic.

\begin{theorem}
Let $f$ is a separable function, i.e., such that there exist factorable functions $f_i: [\underline x_i, \overline x_i] \to \mathbb R$ with
\[
f(x) = \sum_{i=1}^n f_i(x_i) \; .
\]
If the remainder bound of all univariate functions in the atom library $\mathcal L$ satisfies $r_g(A) = 0$ whenever the input model is separable (this condition is satisfied for all operations in Table~\ref{tab::Bounds}), then the interval superposition model $F_{f,X}$ computed by the above outlined arithmetic rules satisfies
\[
\max_{x \in X} \, \left\{  \diam{F_{f,X}(x)} \right\} \; \leq \; \order{ \frac{1}{N}}
\]
for all bounded domains $X \subseteq \mathbb I^n$.
\end{theorem}

\proof
Since the functions $f_i$ depend on one variable only, all intermediate models remain separable (see Lemmas~\ref{lem::zeroRemainder} and~\ref{lem::zeroRemainder}), i.e., we have $r_g(A) = 0$ during the whole evaluation. 

\hfill \proofbox

\section{Implementation and Examples}
\label{sec::examples}
The goal of this section is to illustrate the potential of the proposed interval
 superposition arithmetic for bounding factorable functions. For this aim, the
proposed interval superposition arithmetic has been implemented in the programming 
language \texttt{Julia}. In order to measure the quality of the proposed
arithmetic, we use the following notation for the Hausdorff distance of a
function $f: X \to \mathbb R^m$ and its enclosure function $F$,
\begin{eqnarray}
\label{eq::overestimation}
d_\mathrm{H}( f(X), F(X) ) = \max_{y \in F(x),} \, \min_{x \in f(X)} \, \Vert x - y \Vert_\infty \; .
\end{eqnarray}
Here, $f(X) = \{ f(x) \mid x \in X \}$ denotes the exact image set of $f$ on $X$
and $\Vert \cdot \Vert_\infty$ denotes the standard $\infty$-norm in $\mathbb R^{n}$.

\subsection{Interval superposition models versus Taylor models}
\label{sec::TMvsISA1}

\begin{figure}
	\hspace{-0.4cm}
	\begin{tabular}{cc}
		\includegraphics[scale=0.35]{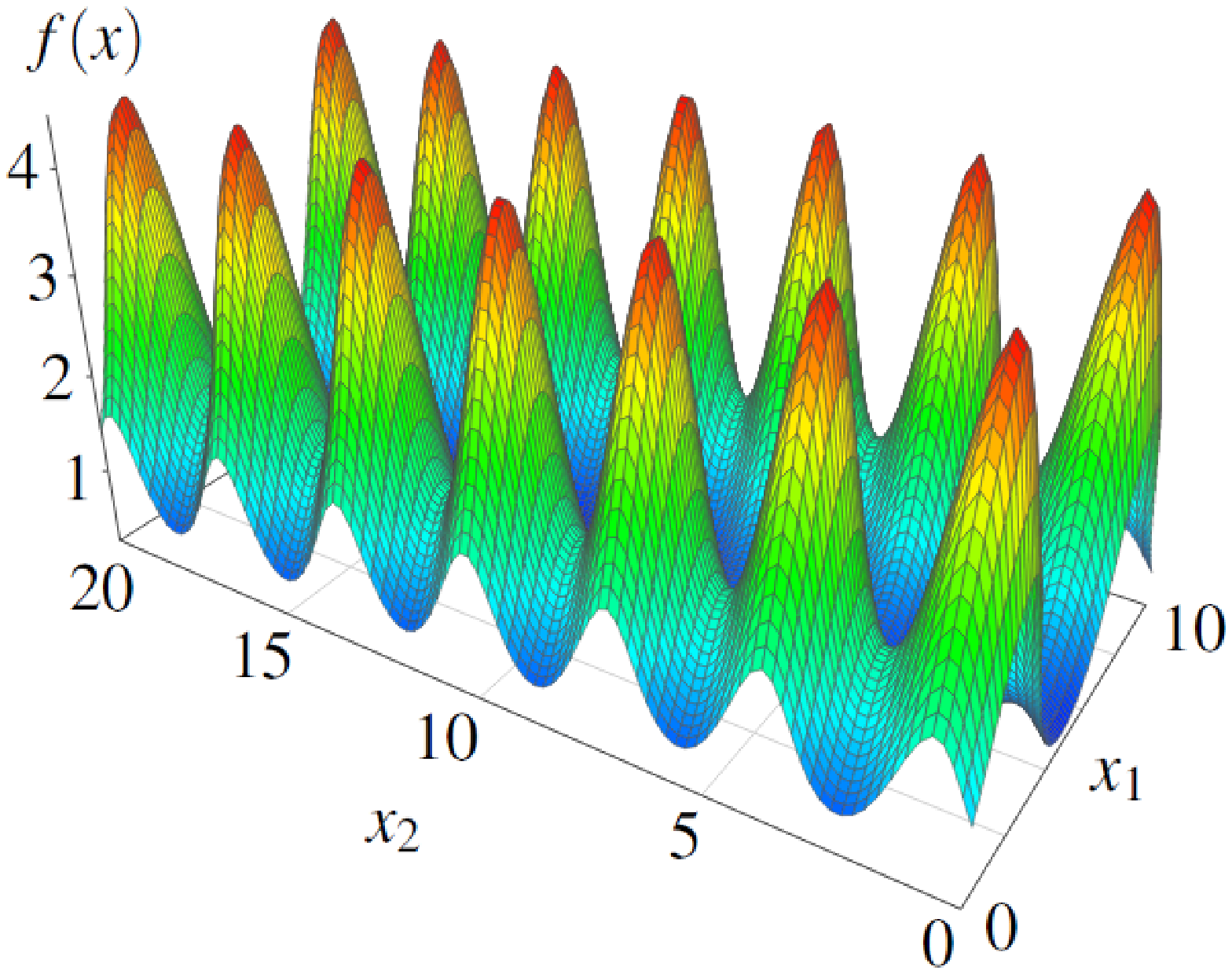} &
		\includegraphics[scale=0.74]{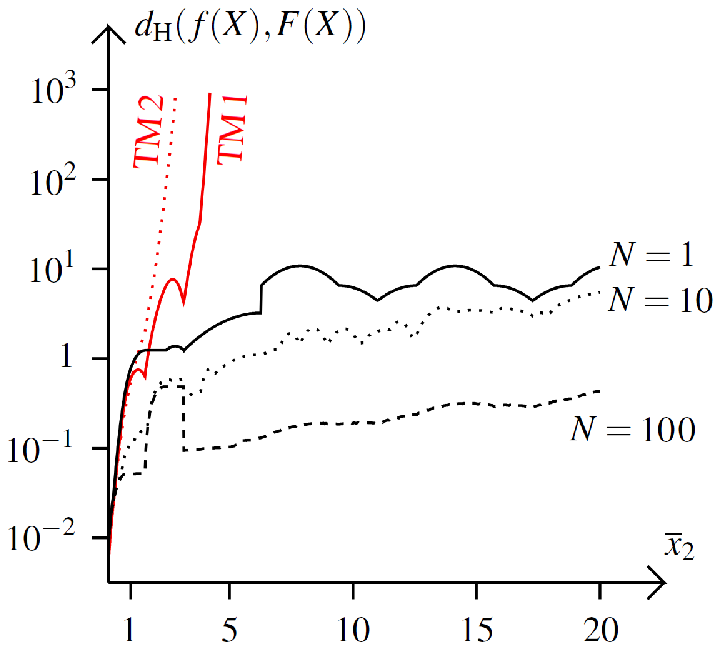} \\
		\includegraphics[scale=0.74]{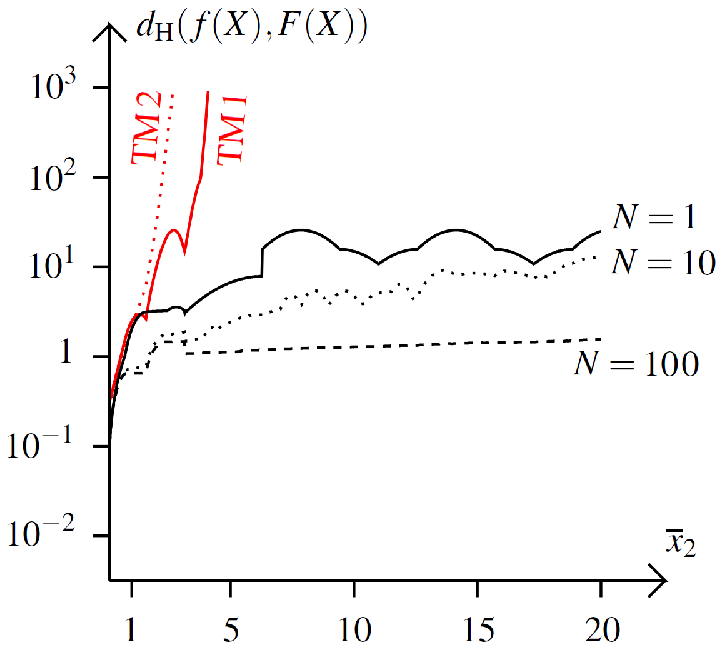} &
		\includegraphics[scale=0.74]{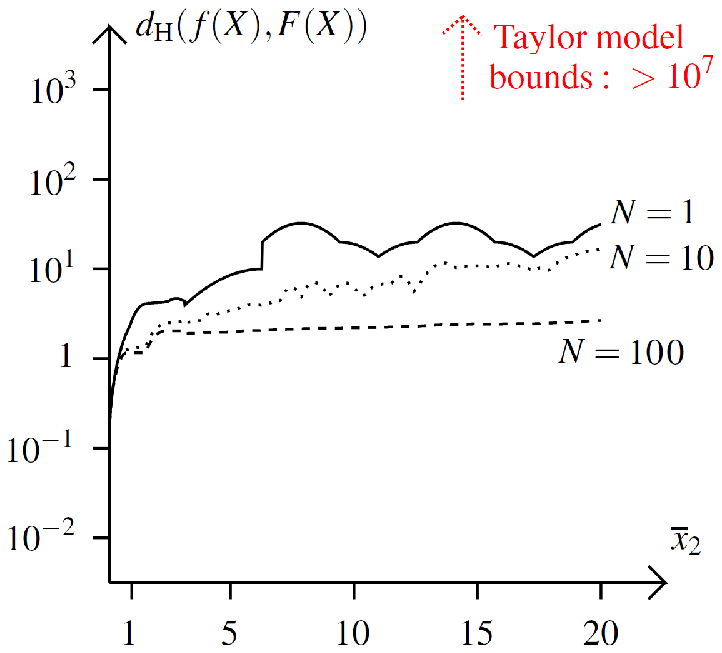}
	\end{tabular}
	\caption{\label{fig::result1}\texttt{Upper-left}: $3$-dimensional visualization of the function $f$ on the domain $X = [0,10] \times [0,20]$. \texttt{Upper-right}: Hausdorff distance between the exact image set $f(X)$ and its enclosure sets on the interval $X = [0,0.1] \times [0,\overline x_2]$ as a function of $\overline x_2 \in [0.1,20]$. \texttt{Lower-left}: Hausdorff distance between the exact image set $f(X)$ and its enclosure sets on the interval $X = [0,1] \times [0,\overline x_2]$ as a function of $\overline x_2 \in [0.1,20]$. \texttt{Lower-right}: Hausdorff distance between the exact image set $f(X)$ and its enclosure sets on the interval $X = [0,10] \times [0,\overline x_2]$ as a function of $\overline x_2 \in [0.1,20]$. In all plots the black solid line corresponds to the results obtained with interval superposition arithmetics with $N=1$. The black dotted lines correspond to interval superposition arithmetic with $N=10$, and the black dashed lines use $N=100$. The red solid and red dotted line correspond to the results obtained with Taylor models of order $1$ and $2$, respectively.}
\end{figure}

The goal of this section is to compare the performance of interval superposition models versus Taylor models on wider domains. Let $f: \mathbb R^2 \to \mathbb R$ denote a non-convex factorable function of the form
\[
f(x) = \exp \left( \sin(x_1) + \sin(x_2) \cos(x_2) \right)
\]
on the two-dimensional domain $X = [0, \overline x_1 ] \times [0, \overline x_2] \subseteq \mathbb R^2$. Here, $\overline x_1 \geq 0$ and $\overline x_2 \geq 0$ are parameters that can be used to control the diameter of the domain $X$. The upper left plot in Figure~\ref{fig::result1} shows a $3$-dimensional visualization of the function $f$ on the interval domain $[0,10] \times [0,20]$, i.e., for $\overline x_1 = 10$ and $\overline x_2 = 20$. The upper right plot in Figure~\ref{fig::result1} shows the overestimation of five different enclosure methods for bounding $f$ for $\overline x_1 = 10^{-1}$ as a function of the domain parameter $\overline x_2 \in [0.1,20]$: the red solid and red dotted lines show the overestimation of Taylor models of order $1$ and $2$, respectively. The black solid, black dotted, and black dashed lines correspond to the overestimation of the enclosures that are obtained by using interval superposition models with $N=1$, $N=10$, and $N=100$. Here, the overestimation is measured in terms of the Hausdorff distance~\eqref{eq::overestimation} between the exact image set and the five different enclosure sets. The lower left plot in Figure~\ref{fig::result1} shows the overestimation of the five mentioned methods for a fixed $\overline x_1 = 1$ as a function of $\overline x_2 \in [0.1,20]$.
Similarly, the lower right plot in Figure~\ref{fig::result1} depicts the corresponding results for $\overline x_1 = 10$, again as a function of $\overline x_2 \in [0.1,20]$. Here, the results for the Taylor models is not shown, as the overestimation error is larger that $10^7$, i.e., Taylor models do not yield reasonable enclosures on this rather large domain. In order to avoid misunderstanding at this point, notice that the width of the exact image set $f(X)$ is monotonically increasing with respect to the parameter $\overline x_2$. However, the Hausdorff difference between $f(X)$ and an enclosure $F(X)$ is not necessarily monotonous in $\overline x_2$. In fact, also the overestimation of standard Taylor models decreases in sections, if the domain $X$ is increased, although one might argue that, overall, a rough trend is that the overestimation error of the enclosure methods increases when increasing the domain~$X$. One aspect that is not shown in the Figure~\ref{fig::result1} is that Taylor models do outperform interval superposition models on very small domains, i.e., if we would zoom in and analyze the overestimation for $\overline x_1, \overline x_2 \leq 10^{-1}$, we could see that Taylor models are the better choice on such small domains. Notice that on the domain 
$[0,10] \times [0,20]$ the overestimation of the interval 
superposition method with $N=100$ yields an enclosure that is approximately $1.62$ 
times larger than the width of the exact range, i.e., the relative over-approximation is approximately $62 \; \%$.
This is in contrast to Taylor models, which yield bounds that are more than $10^7$ times larger than the exact image set. The performance of Taylor models of order larger than $2$ is not shown in the figure, as they perform even worse than the Taylor models of order $1$ and $2$ on the analyzed, particularly large domains. Here, of course, if we would zoom in on smaller domains $X$, we could see that increasing the Taylor model order does improve the accuracy for such smaller $X$~\cite{Berz1997,Berz1998,Sahlodin2011}. 

In order to illustrate how the proposed interval superposition arithmetics performs for another, more challenging example, we introduce the function
\begin{eqnarray}
f_1(x) = \left(
\begin{array}{c}
p_1 \left( e^{-\sin(4x_1)+x_2-x_2^2-x_1^2} - 1 \right) \\[0.16cm]
p_1 \cos \left( \frac{1}{p_1} x_2+p_2 \tan(p_2 x_3) \right) - 2 p_2 x_2^2 \\[0.16cm]
p_1^2 \sin(\cos(x_3))
\end{array}
\right) \; .
\label{eq::f1}
\end{eqnarray}
Notice that $f_1: \mathbb R^{3} \to \mathbb R^3$ is a multivariate non-convex function with parameters $p_1 = \frac{1}{10}$ and $p_2 = \frac{1}{5}$.
In the next step we define the functions 
\begin{eqnarray}
\forall k \in \mathbb N, \qquad f_{k+1}(x) = f_1(f_{k}(x))
\label{eq::fk}
\end{eqnarray}
recursively. The goal of this section is to find enclosure sets of the exact image sets $f_k(X)$ on the rather large interval domain $X = [-0.25\pi,0.25\pi] \times  [-0.5\pi,0.5\pi]^2 \subseteq \mathbb R^3$. The exact image set of the above recursion satisfies a convergence rate condition of the form
\[
\lim_{k \to \infty} \diam{ f_k(X) } \; = \; 0 \; ,
\]
i.e., the diameter of the exact image set contracts to $0$ for sufficiently large $k$.
\begin{figure}
	\begin{center}
		\includegraphics[scale=0.75]{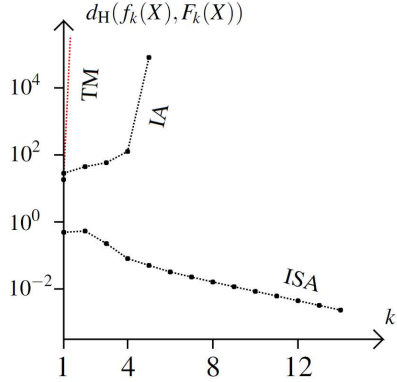}
	\end{center}
	\caption{\label{fig::recursion} The Hausdorff distance between the exact image sets $f_k(X)$ and their computed enclosure sets $F_k(X)$ in dependence on the running index $k$. The corresponding results for standard interval arithmetic are labeled as ``IA''. The other results are obtained by using interval superposition models with $N=20$ labeled as ``ISA'' and first order Taylor models with interval remainder as ``TM'', respectively.}
\end{figure}
Figure~\ref{fig::recursion} shows the Hausdorff distance between the exact image set $f_k(X)$ and the enclosure sets that are obtained by applying standard interval arithmetic, interval superposition arithmetic with $N = 20$, and first order Taylor models with interval remainder bounds. All results are shown in dependence on $k$. Notice that the interval superposition arithmetic yields convergent enclosure sets that are much less conservative than the enclosures that are obtained by Taylor models and standard interval arithmetic. Taylor models of higher expansion orders perform even worse on this example and are therefore not shown in the figure.

\section{Conclusions}
\label{sec::conclusions}
This paper has introduced interval superposition arithmetic and illustrated its 
advantages compared to existing enclosure methods for factorable functions on 
wider domains. The construction of interval superposition models is based on 
derivative-free composition rules which exploit global algebraic properties 
of factorable functions. Interval superposition arithmetic has polynomial 
run-time and storage complexity of order $\order{n N}$, which depends on the 
number $n$ of variables of the factorable functions and the branching accuracy 
$N$. Moreover, this paper has established local and global convergence estimates
of the proposed arithmetic.
From a practical perspective, the main advantage of interval superposition 
arithmetics compared to other enclosure methods is that it yields reasonably 
accurate bounds of the image set of factorable functions on wider interval 
domains, for which existing methods often yield divergent or very conservative 
bounds. This advantage has been illustrated through numerical case studies.

\appendix
\section{Derivation of the remainder bounds from Table~\ref{tab::Bounds}}
\label{sec::bounds}

This section briefly discusses how to derive remainder bounds for
interval superposition arithmetic. These remainder bounds are needed in 
Algorithm~1 and are required to satisfy
\begin{eqnarray}
\left| \sum_{i=1}^n g( \omega + \delta_i ) - (n-1)g(\omega) - g \left( \omega + \sum_{i=1}^n \delta_i \right) \right| \; \leq \; r_g(A)
\end{eqnarray}
for all $\delta \in \mathbb R^n$ with $\forall i \in \{ 1, \ldots, n \}, \; \; L(A_i) \leq a_i + \delta_i \in U(A_i)$. Recall that $g \in \mathcal L$ denotes a univariate atom function and its associated remainder bound $r_g$ depends on the particular properties of $g$. Also recall that we use shorthand $\omega = \sum_{i=1}^n a_i$ as introduced in the first step of Algorithm~1.

\subsection{Exponential}
For the atom function $g(x) = e^x$ we have to bound the expression
\begin{eqnarray}
\sum_{i=1}^n g( \omega + \delta_i ) - (n-1)g(\omega) - g \left( \omega + \sum_{i=1}^n \delta_i \right) &=& e^{\omega} \left[ \sum_{i=1}^n e^{\delta_i} - (n-1) - \prod_{i=1}^n e^{\delta_i} \right] \notag
\end{eqnarray}
for all $\delta_i$ with $L(A_i) \leq a_i + \delta_i \leq U(A_i)$. Let us apply the addition theorem for the exponential function,
$$e^{\omega + \delta_i} = e^{\omega} e^{\delta_i} \quad \text{and} \quad e^{\omega + \sum_{i=1}^n \delta_i} = e^{\omega} \prod_{i=1}^n e^{\delta_i} \; .$$
It is convenient to introduce the auxiliary variables $t_i = e^{\delta_i} - 1$ such that
\begin{eqnarray}
e^{\omega} \left[ \sum_{i=1}^n e^{\delta_i} - (n-1) - \prod_{i=1}^n e^{\delta_i} \right] &=& e^{\omega} \left[ \sum_{i=1}^n t_i + 1 - \prod_{i=1}^n (1+t_i) \right] \; .
\end{eqnarray}
The absolute value of this expression can be bounded as
\[
e^{\omega} \left| \sum_{i=1}^n t_i + 1 - \prod_{i=1}^n (1+t_i) \right| \; \leq \; e^{\omega} \left( \prod_{i=1}^n (1+s_i) - \sum_{i=1}^n s_i - 1 \right)  \; 
\]
with $s_i = \max \left\{ e^{\mathrm{U}\left( A_i \right)-a_i} - 1, 1 - e^{\mathrm{L}\left( A_i \right)-a_i}  \right\}$. This motivates to choose the central points $a_i = \log \left( \frac{1}{2} \left( e^{\mathrm{U}\left( A_i \right)} + e^{\mathrm{L}\left( A_i \right)} \right) \right)$ such that $s_i$ takes the smallest possible value, given by
\[
s_i = \frac{e^{\mathrm{U}(A_i)}-e^{\mathrm{L}\left( A_i \right)}}{e^{\mathrm{U}\left( A_i \right)} + e^{\mathrm{L}\left( A_i \right)}} \; .
\]
In summary, we have shown that
\begin{eqnarray}
\left| \sum_{i=1}^n g( \omega + \delta_i ) - (n-1)g(\omega) - g \left( \omega + \sum_{i=1}^n \delta_i \right) \right| &\leq&  e^{\omega} \left( \prod_{i=1}^n (1+s_i) - \sum_{i=1}^n s_i - 1 \right) \; = \; r_g(A) \; . \notag
\end{eqnarray}

\subsection{Inverse}
The aim of this section is to find a remainder bound for the atom function 
$g(x) = \frac{1}{x}$ on the positive domain $\mathbb R_{++} = \{ x \mid x > 0 \}$.
Bounds on the domain $\mathbb R_{--} = \{ x \mid x < 0 \}$  can be found 
analogously. If an interval contains $0$, the bounds are set to $[-\infty,\infty]$.
We start with the equation
\begin{align*}
&\sum_{i=1}^n g( \omega + \delta_i ) - (n-1)g(\omega) - g \left( \omega + \sum_{i=1}^n \delta_i \right) 
= \sum_{i=1}^n \frac{1}{\omega+\delta_i} - \frac{1}{\omega + \sum_{i=1}^n \delta_{i}} - \frac{n-1}{\omega} \\ 
&\quad = \frac{1}{\omega} \left( \sum_{i=1}^n \frac{-\delta_i}{\omega+\delta_i} + \frac{\sum_{i=1}^n \delta_i}{\omega + \sum_{i=1}^n \delta_{i}} \right) 
= \frac{1}{\omega} \, \frac{1}{\omega + \sum_{i=1}^n \delta_{i}} \, \left( \sum_{i=1}^n \frac{ \delta_i ( \delta_i - \sum_{k=1}^n \delta_k)}{\omega+\delta_i} \right) \; . 
\end{align*}
Next, we bound the terms in the last equation separately under the assumption that 
$\lambda(A) > 0$,
\begin{align*}
\left| \frac{1}{\omega + \sum_{i=1}^n \delta_{i}} \right| \leq \frac{1}{\lambda(A)} \; , \quad 
\left| \frac{\delta_i}{\omega + \delta_{i}} \right| &\leq \max \left\{ \, \frac{a_i - L(A_i)}{ \omega - a_i + L(A_i) } \, , \, \frac{U(A_i) - a_i}{ \omega - a_i + U(A_i)} \, \right\} \; = \;  s_i \; , \\
\text{and} \quad  \left| \delta_i - \sum_{k=1}^n \delta_k \right| &\leq \mu(A) - \omega - (U(A_i)-a_i) \; . 
\end{align*}
Substituting these inequalities yields the desired remainder bound
\[
\left| \sum_{i=1}^n g( \omega + \delta_i ) - (n-1)g(\omega) - g \left( \!\omega + \sum_{i=1}^n \delta_i \!\right) \right|  \leq   \frac{\sum_{i=1}^n s_i \left( \mu(A) - \omega - (U(A_i)-a_i) \right)}{\omega \lambda(A)}  =  r_{g}(A). 
\]

\subsection{Logarithm}
The aim of this section is to find a remainder bound for the atom function $g(x) = \log(x)$ on the positive domain $\mathbb R_{++} = \{ x \mid x > 0 \}$,
\begin{align*}
&\sum_{i=1}^n g( \omega + \delta_i ) - (n-1)g(\omega) 
- g \left( \omega + \sum_{i=1}^n \delta_i \right) \\
& \qquad = \sum_{i=1}^n \log \left( \omega+\delta_i \right) 
- \log \left( \omega + \sum_{i=1}^n \delta_{i} \right) 
- (n-1)\log \left( \omega \right) \notag \\
&\qquad = \log \left( \frac{\prod_{i=1}^n 
\left( \omega+\delta_i \right) }{\omega^{n-1} \left( \omega + \sum_{i=1}^n \delta_{i} \right)} \right) = 
\log \left( 1 + \frac{\prod_{i=1}^n \left( \omega+\delta_i \right) - \omega^{n-1} \left( \omega + \sum_{i=1}^n \delta_{i} \right) }{\omega^{n-1} \left( \omega + \sum_{i=1}^n \delta_{i} \right)} \right) \; . \notag
\end{align*}
The desired bound is found by bounding the absolute value of this term, choosing
the central points 
$a_i = \frac{U(A_i)+L(A_i)}{2}$ such that $\left| \delta_i \right| \leq s_i = \frac{U(A_i) - L(A_i)}{2}$ and
\begin{align*}
&\left| \sum_{i=1}^n g( \omega + \delta_i ) - (n-1)g(\omega) - g \left( \omega + \sum_{i=1}^n \delta_i \right) \right| \\
&\qquad\leq - \log \left( 1 - \frac{\prod_{i=1}^n \left( \omega+s_i \right) - \omega^{n-1} \left( \omega + \sum_{i=1}^n s_{i} \right) }{\omega^{n-1} \lambda(A)} \right) = r_g(A) \; . \notag
\end{align*}

\subsection{Sine and Cosine}
In order to derive remainder bounds for the sine and cosine functions
we use Euler's formula, $e^{\imag x} = \cos(x) + \imag \sin(x)$ with 
$\imag = \sqrt{-1}$. The derivation requires the following steps.

\bigskip
\noindent
\textbf{Step 1.}
In the first step, we derive for all $k \in \{ 1, \ldots, n \}$ the bound
\begin{align*}
\left| e^{\pm \imag \delta_k} - 1 \right| &= \left| \cos(\pm \delta_k)-1 + \imag \sin( \pm \delta_k) \right| \notag \\
&= 2 \left| \sin \left( \pm \frac{\delta_k}{2} \right) \right| 
\leq 2 \left| \sin \left( \left[ -\frac{U(A_k)-L(A_k)}{4},\frac{U(A_k)-L(A_k)}{4} \right] \right) \right| = s_k \; .
\end{align*}
Here, the expression for the scalars $s_k$ is evaluated by using standard interval arithmetic, i.e., 
\begin{align*}
s_k &= 2 \left| \sin \left( \left[ -\frac{U(A_k)-L(A_k)}{4},\frac{U(A_k)-L(A_k)}{4} \right] \right) \right| \\ &= 
\begin{cases}
2 \sin \left( \frac{U(A_k)-L(A_k)}{4} \right)  &\text{if} \ \
\frac{U(A_k)-L(A_k)}{4} \leq \frac{\pi}{2} \\
2 &\text{otherwise} \!\!\!
\end{cases}\;.
\end{align*}

\bigskip
\noindent
\textbf{Step 2.}
In the second step, we use the bounds $s_k$ to derive the auxiliary inequalities
\begin{eqnarray}
\left| \sum_{k=1}^n e^{\pm \imag \delta_k} - \prod_{k=1}^n e^{\pm \imag \delta_k} - (n-1) \right| &\leq& \prod_{k=1}^n (1+s_k) - \sum_{k=1}^n s_k - 1 \; . \notag 
\end{eqnarray}

\bigskip
\noindent
\textbf{Step 3.}
The auxiliary inequalities from Step 2 are used to establish the inequalities
\begin{align*}
&\left| \sum_{k=1}^n \cos \left( \delta_k \right) - \cos \left( \sum_{k=1}^n \delta_{k} \right) - (n-1) \right| \\
& \ = \frac{1}{2} \left| \sum_{k=1}^n e^{\imag \delta_k} + \sum_{k=1}^n e^{-\imag \delta_k} - \prod_{k=1}^n e^{\imag \delta_k} - \prod_{k=1}^n e^{-\imag \delta_k} - 2(n-1)  \right| 
\leq \prod_{k=1}^n (1+s_k) - \sum_{k=1}^n s_k - 1 
\end{align*}
and, using an analogous argument,
\begin{eqnarray}
\left| \sum_{k=1}^n \sin \left( \delta_k \right) - \sin \left( \sum_{k=1}^n \delta_{k} \right) \right| \; \leq \; \prod_{k=1}^n (1+s_k) - \sum_{k=1}^n s_k - 1 \;. \notag 
\end{eqnarray}

\bigskip
\noindent
\textbf{Step 4.}
For the sine function, the estimate from Step~3 yields the remainder bound
\begin{align*}
R_{\sin}(\delta) &= \left| \sum_{k=1}^n g( \omega + \delta_k ) - (n-1)g(\omega) - g \left( \omega + \sum_{k=1}^n \delta_k \right) \right|  \\
&= \left| \sum_{k=1}^n \sin \left( \omega+\delta_k \right) - \sin \left( \omega + \sum_{k=1}^n \delta_{k} \right) - (n-1) \sin \left( \omega \right) \right| \notag \\
&= \left| \hphantom{{}+{}} \sin(\omega) \left( \sum_{k=1}^n \cos \left( \delta_k \right) - \cos \left( \sum_{k=1}^n \delta_{k} \right) - (n-1) \right)  \right . \\
&\hphantom{{}={}} \, \left . + \cos(\omega) \left( \sum_{k=1}^n \sin \left( \delta_k \right) - \sin \left( \sum_{k=1}^n \delta_{k} \right) \right) \right|  \\
&\leq \left( |\sin(\omega)| + |\cos(\omega)| \right) \left( \prod_{k=1}^n (1+s_k) - \sum_{k=1}^n s_k - 1 \right) = r_g(A) \; . 
\end{align*}
Similarly, the corresponding bound for the cosine function is given by
\begin{eqnarray}
R_{\cos}(\delta) &=& \left| \sum_{k=1}^n g( \omega + \delta_k ) - (n-1)g(\omega) - g \left( \omega + \sum_{k=1}^n \delta_k \right) \right| \notag \\
&=& \left| \sum_{k=1}^n \cos \left( \omega+\delta_k \right) - \cos \left( \omega + \sum_{k=1}^n \delta_{k} \right) - (n-1) \cos \left( \omega \right) \right| \notag \\
&=& \left| \cos(\omega) \left( \sum_{k=1}^n \cos \left( \delta_k \right) - \cos \left( \sum_{k=1}^n \delta_{k} \right) - (n-1) \right) - \sin(\omega) \left( \sum_{k=1}^n \sin \left( \delta_k \right) - \sin \left( \sum_{k=1}^n \delta_{k} \right) \right) \right| \notag \\
&\leq& \left( |\sin(\omega)| + |\cos(\omega)| \right) \left( \prod_{k=1}^n (1+s_k) - \sum_{k=1}^n s_k - 1 \right) = r_g(A) \; . \notag
\end{eqnarray}

\subsection{Tangent}
In order to construct a remainder bound for the function $g(x) = \tan(x)$ on the
open domain $\left( -\frac{\pi}{2}, \frac{\pi}{2} \right)$ it is helpful to 
notice that the addition theorem for this function, 
\[
\tan(x + y) = \frac{\tan(x) + \tan(y)}{1 - \tan(x)\tan(y)} \; ,
\]
can alternatively be written in the difference form
\begin{eqnarray}
\label{eq::tangentDifference}
\tan(x + y) - \tan(x) - \tan(y) &=& \tan(x) \tan(y) \tan(x+y) \; .
\end{eqnarray}
The correctness of this equation can be verified by multiplying the addition theorem for the tangent function by $1 - \tan(x)\tan(y)$ on both sides and by re-bracketing terms. A generalization of the difference formula~\eqref{eq::tangentDifference} for general sums is given by the equation
\begin{eqnarray}
\label{eq::GenTanDiff}
\rho_0(\delta) \; = \; \tan \left( \sum_{i=1}^n \delta_i \right) - \sum_{i=1}^n \tan( \delta_i ) &=& \sum_{i=1}^{n-1} \tan(\delta_{i+1}) \tan \left( \sum_{k=1}^i \delta_k \right)\tan \left( \sum_{k=1}^{i+1} \delta_k \right);\, 
\end{eqnarray}
which is proven by induction. For $n=2$,~\eqref{eq::GenTanDiff} reduces to~\eqref{eq::tangentDifference}. For the induction step, we have
\begin{align*}
&\tan \left( \sum_{i=1}^{n+1} \delta_i \right) - \sum_{i=1}^{n+1} \tan( \delta_i )\\
&\qquad\qquad= \tan \left( \sum_{i=1}^{n+1} \delta_i \right) - \tan \left( \sum_{i=1}^{n} \delta_i \right)  - \tan(\delta_{n+1}) + \tan \left( \sum_{i=1}^{n} \delta_i \right)  - \sum_{i=1}^{n} \tan( \delta_i ) \notag \\
&\qquad\qquad\overset{\eqref{eq::tangentDifference}}{=} \tan \left( \sum_{i=1}^{n+1} \delta_i \right) \tan \left( \sum_{i=1}^{n} \delta_i \right) \tan(\delta_{n+1}) + \left[ \tan \left( \sum_{i=1}^{n} \delta_i \right)  - \sum_{i=1}^{n} \tan( \delta_i ) \right] \notag \\
&\qquad\qquad\overset{\eqref{eq::GenTanDiff}}{=} \sum_{i=1}^{n} \tan(\delta_{i+1}) \tan \left( \sum_{k=1}^i \delta_k \right)\tan \left( \sum_{k=1}^{i+1} \delta_k \right) \; .
\end{align*}
Thus, the difference formula~\eqref{eq::GenTanDiff} holds for all integers $n$. 
In order to generalize the above formula further for the case $\omega \neq 0$, 
the following algebraic manipulations are made 
\begin{align*}
\rho(\delta) &= g \left( \omega + \sum_{i=1}^n \delta_i \right) + (n-1)g(\omega) - \sum_{i=1}^n g( \omega + \delta_i ) \\
&= \left[ \tan \left( \omega + \sum_{i=1}^n \delta_i \right) - \tan(\omega) \right] - \sum_{i=1}^n \left[ \tan( \omega + \delta_i ) - \tan(\omega) \right]  \\
&\overset{\eqref{eq::tangentDifference}}{=} 
\tan \left( \sum_{i=1}^n \delta_i \right) \left[ 1 + \tan \left( \omega + \sum_{i=1}^n \delta_i \right) \tan(\omega) \right] - \sum_{i=1}^n \tan( \delta_i ) \left[ 1 + \tan( \omega + \delta_i ) \tan(\omega) \right] \\
&= \hphantom{{}+{}} \left( \tan \left( \sum_{i=1}^{n} \delta_i \right) - \sum_{i=1}^{n} \tan( \delta_i ) \right)\\
&\hphantom{{}={}}+ \tan(\omega) \left( \tan \left( \sum_{i=1}^n \delta_i \right)\tan \left( \omega + \sum_{i=1}^n \delta_i \right) - \sum_{i=1}^n \tan( \delta_i ) \tan( \omega + \delta_i ) \right)\\
&\overset{\eqref{eq::GenTanDiff}}{=} \rho_0(\delta) + \tan(\omega) \left( \hphantom{{}+{}} \rho_0(\delta) \tan \left( \omega + \sum_{i=1}^n \delta_i \right) \right .\\
&\hphantom{{}={}} \kern2.7cm + \left . \sum_{i=1}^n \tan \left(  \delta_i \right) \left[ \tan \left( \omega + \sum_{i=1}^n \delta_i \right) - \tan( \omega + \delta_i ) \right] \right) 
 \\
&=\hphantom{{}+{}}\rho_0(\delta) \left[1 + \tan(\omega) \tan \left(\omega + \sum_{i=1}^n \delta_i \right) \right] \\ 
& \hphantom{{}={}}+ \sum_{i=1}^n \tan(\omega) \tan \left(  \delta_i \right)  \tan \left( \sum_{k \neq i} \delta_k \right) \left[ 1 + \tan(\omega+\delta_i)\tan \left( \sum_{k \neq i} \delta_k \right) \tan \left(\omega + \sum_{i=1}^n \delta_i \right) \right] \; .
\notag
\end{align*}
The right-hand expression can be bounded with interval arithmetic yielding
\begin{eqnarray}
\rho(\delta) \leq r_g(A) &=& \left| \sum_{i=1}^{n-1} \tan( S_{i+1} ) \tan \left( \sum_{k=1}^i S_k \right)\tan \left( \sum_{k=1}^{i+1} S_k \right)\left[1 + \tan(\omega) \tan \left( \omega + \Sigma \right) \right] \right. \notag \\[0.16cm]
& & \left. + \sum_{i=1}^n \tan(\omega) \tan \left( S_i \right)  \tan \left( T_i \right) \left[ 1 + \tan( \omega+S_i)\tan \left( T_i \right) \tan \left( \omega + \Sigma \right) \right] \right| \; , \notag
\end{eqnarray}
the desired bound. Here we have introduced the auxiliary variables
\[
s_i = \frac{U(A_i)-L(A_i)}{2} \; , \; \; S_i = [-s_{i},s_{i}]  \quad \text{and} \quad \sigma = \sum_{i=1}^n s_i  \; , \; \; \Sigma = [-\sigma, \sigma]  \; , \; \; T_i = [-\sigma + s_i, \sigma-s_i] \; . 
\]

\section*{Acknowledgments}
\small
This research was supported by
National Natural Science Foundation China (NSFC), Nr.~61473185,
as well as ShanghaiTech University, Grant-Nr.~\mbox{F-0203-14-012}.
\normalsize


\begin{thebibliography}{99}
	
	
	
	\bibitem{Battles2004}
	Z.~Battles, L.N.~Trefethen.
	\newblock An extension of {MATLAB} to continuous functions and operators.
	\newblock SIAM J.~Sci.~Comput.~25:1743--1770, 2004.
	
	\bibitem{Berz1997}
	M.~Berz.
	\newblock From Taylor series to Taylor models.
	\newblock In Nonlinear Problems in Accelerator Physics, American Institute of Physics CP405, pp.:1--27, 1997.
	
	\bibitem{Berz1998}
	M.~Berz, G.~Hoffst{\"a}tter.
	\newblock Computation and application of Taylor polynomials with remainder bounds.
	\newblock Reliab.~Comput.~4:83--97, 1998.
	
	\bibitem{Bompadre2013}
	A.~Bompadre, A.~Mitsos, B.~Chachuat.
	\newblock Convergence analysis of Taylor and McCormick-Taylor models.
	\newblock Journal of Global Optimization 57(1):75--114, 2013.
	
	
	\bibitem{Chachuat2015}
	B.~Chachuat, B.~Houska, R.~Paulen, N.~Peric, J.~Rajyaguru, M.E.~Villanueva.
	\newblock Set theoretic approaches in analysis, estimation and control of nonlinear systems.
	\newblock IFAC-PapersOnLine Volume 48(8), pp:981--995, 2015.
	
%
	
	
	\bibitem{Eckmann1984}
	J.P.~Eckmann, H.~Koch, P.~Wittwer.
	\newblock A computer-assisted proof of universality in area-preserving maps.
	\newblock Memoirs of the AMS 47:289, 1984.
	
	\bibitem{Figueiredo2004}
	L.H.~de~Figueiredo, J.~Stolfi.
	\newblock Affine arithmetic: Concepts and applications.
	\newblock Numerical Algorithms 37(1-4):147--158, 2004.
	
	\bibitem{Floudas2007}
	C.A. Floudas and O.~Stein.
	\newblock The Adaptative Convexification Algorithm: a Feasible Point Method for Semi-Infinite Programming.
	\newblock SIAM Journal on Optimization, 18(4):1187--1208, 2007.
	
	\bibitem{Floudas2013}
	C.A.~Floudas.
	\newblock Deterministic global optimization: theory, methods and applications.
	\newblock Springer Science \& Business Media, Vol. 37, 2013.
	
	\bibitem{Houska2015}
	B.~Houska, M.E.~Villanueva, B.~Chachuat.
	\newblock Stable Set-Valued Integration of Nonlinear Dynamic Systems using Affine Set Parameterizations.
	\newblock SIAM Journal on Numerical Analysis, 53(5), pp:2307--2328, 2015.
	
	\bibitem{Kurzhanski1997}
	A.~Kurzhanski, I.~Valyi.
	\newblock Ellipsoidal Calculus for Estimation and Control.
	\newblock Series in Systems \& Control: Foundations \& Applications, Birkh\"{a}user, 1997.
	
	\bibitem{Lasserre2009}
	J.B. Lasserre.
	\newblock Moments, Positive Polynomials and Their Applications.
	\newblock Imperial College Press, 2009.
	
	\bibitem{Rokne1995}
	Q.~Lin, J.G.~Rokne.
	\newblock Methods for bounding the range of a polynomial.
	\newblock J.~Comput Appl Math 58:193--199, 1995.
	
	\bibitem{Neher2007}
	M.~Neher, K.R.~Jackson, N.S.~Nedialkov. On Taylor model based integration of ODEs.
	\newblock SIAM Journal on Numerical Analysis 45:236--262, 2007.
	
	\bibitem{Neumaier2004}
	A.~Neumaier.
	\newblock Complete search in continuous global optimization and constraint satisfaction.
	\newblock Acta Numer.~13:271--369, 2004.
	
	\bibitem{Makino1999}
	K.~Makino, M.~Berz.
	\newblock Efficient control of the dependency problem based on Taylor model methods.
	\newblock Reliab.~Comput.~5(1):3--12, 1999.
	
	\bibitem{McCormick1976}
	G.P.~McCormick.
	\newblock Computability of global solutions to factorable nonconvex programs: Part I -- Convex underestimating problems. Mathematical Programing 10:147--175, 1976.
	
	\bibitem{Misener2013}
	R.~Misener, C.A.~Floudas.
	\newblock GloMIQO: Global Mixed-Integer Quadratic Optimizer.
	\newblock Journal of Global Optimization, 57(1):3--50, 2013.
	
	\bibitem{Mitsos2008}
	A.~Mitsos, P.~Lemonidis, P.I.~Barton.
	\newblock Global solution of bilevel programs with a nonconvex inner program.
	\newblock Journal of Global Optimization 42.4:475--513, 2008.
	
	\bibitem{Mitsos2009}
	A.~Mitsos, B.~Chachuat, P.I.~Barton.
	\newblock McCormick-based relaxations of algorithms.
	\newblock SIAM Journal on Optimization 20(2):573--601, 2009.
	
	\bibitem{Moore1966}
	R.E.~Moore.
	\newblock Interval Analysis.
	\newblock Prentice-Hall, Englewood Cliffs, NJ, 2966.
	
	\bibitem{Moore2009}
	R.E.~Moore, R.B.~Kearfott, M.J.~Cloud.
	\newblock Introduction to Interval Analysis.
	\newblock SIAM, Philadelphia, PA, 2009.
	
	\bibitem{Rajyaguru2016}
	J.~Rajyaguru, M.E.~Villanueva, B.~Houksa, B.~Chachuat.
	\newblock Higher-Order Inclusions of Factorable Functions by Chebyshev Models.
	\newblock Journal of Global Optimization, Volume~68(2), pp.~413--438, 2017.
	
	\bibitem{Ratschek1984}
	H.~Ratschek, J.~Rokne.
	\newblock Computer Methods for the Range of Functions.
	\newblock Series in Mathematics and Its Applications, Ellis Horwood Ltd, Mathematics and Its Applications, Chichester, UK, 1984.
	
	\bibitem{Sahinidis1996}
	N.V.~Sahinidis.
	\newblock A general purpose global optimization software package.
	\newblock Journal of Global Optimization, 8(2):201--205, 1996.
	
	\bibitem{Sahlodin2011}
	A.M.~Sahlodin, B.~Chachuat.
	\newblock Convex/concave relaxations of parametric ODEs using Taylor models.
	\newblock Computers and Chemical Engineering 35(5):844--857, 2011.
	
	\bibitem{Tawarmalani2005}
	M.~Tawarmalani, N.V.~Sahinidis.
	\newblock A polyhedral branch-and-cut approach to global optimization.
	\newblock Mathematical Programming, 103(2):225--249, 2005.
	
	\bibitem{Trefethen2007}
	L.N.~Trefethen.
	\newblock Computing numerically with functions instead of numbers.
	\newblock Math.~Comput.~Sci.~1:9--19, 2007.
	
	\bibitem{Townsend2013}
	A.~Townsend, L.N.~Trefethen.
	\newblock An extension of {C}hebfun to two dimensions.
	\newblock SIAM J.~Sci.~Comput 35(6):C495--C498, 2013.
	
	\bibitem{Villanueva2015}
	M.E.~Villanueva, J.~Rajyaguru, B.~Houska, B.~Chachuat.
	\newblock Ellipsoidal arithmetic for multivariate systems.
	\newblock Comput.~Aided Chem.~Eng.~37:767--772, 2015.
	
	\bibitem{Villanueva2015a}
	M.E.~Villanueva, B.~Houska, B.~Chachuat.
	\newblock Unified Framework for the Propagation of Continuous-Time Enclosures for Parametric Nonlinear ODEs.
	\newblock J.~of Global Optim 62(3), pp:575--613, 2015.
	
\end{thebibliography}

\end{document}